\documentclass[a4paper]{article}
\usepackage{amsmath,amssymb,amsfonts}

\newtheorem{theorem}{Theorem}[section]

\title{On the G\"odel's formula}
\author{Jailton C. Ferreira}

\date{ }
\begin{document}
\maketitle \pagenumbering{arabic}

\begin{abstract}
The proof of G\"odel's  first incompleteness theorem includes the
construction of an arithmetic formula $G$ that represents the
metamathematical statement: the formula $G$ is not provable. This
article examines the formula $G$ (of G\"odel). We demonstrated
that there is no Gödel's number for the formula $G$ if number of provable well formed
formulae with one free variable is finite. If there is a non-finite number of provable propositions in the theory, then Gödel's formula also does not possess the Gödel's number.
\end{abstract}

\section{Introduction} \label{sec-1}

\hspace{22pt} The G\"odel's formula is frequently comprehended as
a self-referential statement like
\begin{center}
\textit{This sentence is not provable. \\
Let $G$ be the name of the above sentence.}
\end{center}

\hspace{22pt} In section \ref{sec2}, two ways of deriving Gödel's formula are given. In the first way, the arguments of substitution operator $sub$ are Gödel's numbers; in the second the arguments of the substitution operator $Sb$ are formulas and a free variable.

\hspace{22pt} In section \ref{sec3} starting from the diagonal
theorem we conclude that if there is an infinite set $\mathbb{S}$
of well formed formulae with one free variable such that the
elements of $\mathbb{S}$ have proofs in the theory $T$ 
then the G\"odel's sentence does not have Gödel's number.

\hspace{22pt} Section \ref{sec4} presents the conclusions.

%
%

\section{G\"odel's formula} \label{sec2}

%
%

\subsection{G\"odel's formula (1)} \label{sec2.1}

\hspace{22pt} A category of names of formulae, denominated
\textit{structural-descriptive names} by Tarski ~\cite{Tarski}, is
applied to names that describe the words that compose the denoted
expression.

\hspace{22pt}  G\"odel's numbering form attributes a distinctive
numeral to each symbol of the alphabet of a formal language. It
possesses an effective method to map each symbol, sequence of
symbols (which can be a well formed formula) or sequence of well
formed formulae (which can be proof of a theorem) in a numeral
(denominated G\"odel's number), and it possesses an effective
method to map each G\"odel's number in the symbol or sequences of
symbols corresponding to the G\"odel's number.

\hspace{22pt} It is evident that the G\"odel's number of a formula
is a structural-descriptive name of the formula. The formula named
$G$ has a second name that is its G\"odel's number. We will now
build the formula $G$.

\hspace{22pt} Let $y$ be the Gödel number of a well-formed formula
with a single occurrence of the only free variable, $z$

\begin{equation}\label{eq2.1-1}
\ulcorner y\urcorner \hspace{44pt} y
\end{equation}
or
\begin{equation}\label{eq2.1-2}
\ulcorner y\urcorner \hspace{44pt} \Delta z \Omega
\end{equation}

where $\Delta$ and $\Omega$ are sequences of symbols.

Let us \textbf{define} the operator $sub(\ulcorner y \urcorner, \ulcorner z \urcorner, \ulcorner j \urcorner)$ as being Gödel's number of the expanded Gödel's number
obtained with the substitution, in the expanded G\"odel's number $\ulcorner y\urcorner$, of the only free variable, $\ulcorner z\urcorner$, by the expanded Gödel's number of $\ulcorner j\urcorner$.

\hspace{22pt} Let us consider the formula
\begin{equation}\label{eq2.1-3}
\ulcorner n \urcorner \hspace{44pt} \neg (\exists \ulcorner r \urcorner: \exists \ulcorner s \urcorner: (P(\ulcorner r \urcorner, \ulcorner s \urcorner) \land (\ulcorner s \urcorner =  sub(\ulcorner y \urcorner, \ulcorner z \urcorner, \ulcorner y \urcorner) ))
\end{equation}

$P(\ulcorner r \urcorner, \ulcorner s \urcorner)$ is true if the sequence of symbols with G\"odel's number
$\ulcorner r \urcorner$ proves the formula with G\"odel's number $\ulcorner s \urcorner$. In English
\eqref{eq2.1-3} means: \textit{there is not a proof  for the formula
whose G\"odel's numeral is $sub(\ulcorner y \urcorner, \ulcorner z \urcorner, \ulcorner y \urcorner)$}.

\hspace{22pt} Assuming the variable $\ulcorner y \urcorner$ of the operator $sub(\ulcorner y \urcorner, \ulcorner z \urcorner, \ulcorner y \urcorner)$ the numerical value $\ulcorner n \urcorner$, the formula in
\eqref{eq2.1-3} changes to

\begin{equation}\label{eq2.1-4}
\ulcorner G \urcorner \hspace{44pt} \neg (\exists \ulcorner r \urcorner: \exists \ulcorner s \urcorner: (P(\ulcorner r \urcorner, \ulcorner s \urcorner) \land (\ulcorner s \urcorner = sub(\ulcorner n \urcorner, \ulcorner z \urcorner, \ulcorner n \urcorner) ))
\end{equation}
The reason for \eqref{eq2.1-4} was to replace a variable with a numeral.

\hspace{22pt} The proof of G\"odel's first incompleteness theorem includes the
construction of an arithmetic formula $G$ that would represent the
metamathematical statement: formula $G$ is not provable. This
formula $G$ is in Hofstadter ~\cite{Hofstadter} or Nagel and
Newman ~\cite{Nagel&Newman}, for instance,

\begin{equation}\label{eq2.1-5}
\neg (\exists \ulcorner r \urcorner: \exists \ulcorner s \urcorner: (P(\ulcorner r \urcorner, \ulcorner s \urcorner) \land (\ulcorner s \urcorner = sub(\ulcorner n \urcorner, \ulcorner z \urcorner, \ulcorner n \urcorner))))
\end{equation}

%
%

\subsection{G\"odel's formula (2)} \label{sec2.2}

\hspace{22pt} Let $\ulcorner y \urcorner$ be the G\"odel's number of a well formed
formula with a single free variable, $z$

\begin{equation}\label{eq2.2-1}
\ulcorner y \urcorner \hspace{44pt} y
\end{equation}

where $y$ is the formula.

\hspace{22pt} Let us define the function $Sb(y, z, j)$ as being the \textbf{formula}
obtained with the substitution in the formula $y$ the only free variable, $z$, by the formula $j$. 

\hspace{22pt}  Let us consider the formula
\begin{equation}\label{eq2.2-2}
\neg (\exists \ulcorner r \urcorner: \exists \ulcorner s \urcorner: (P(\ulcorner r \urcorner, \ulcorner s \urcorner) \land (\ulcorner s \urcorner = \ulcorner Sb(y,z, y) \urcorner ))
\end{equation}

$P(\ulcorner r \urcorner, \ulcorner s \urcorner)$ is true if the sequence of symbols with G\"odel's number
$\ulcorner r \urcorner$ proves the formula with G\"odel's number $\ulcorner s \urcorner$. In English
\eqref{eq2.2-2} means: \textit{there is not a proof  for the formula
whose G\"odel's numeral is $ \ulcorner Sb(y, z, y) \urcorner $ }.

\hspace{22pt} Let $\ulcorner n \urcorner$ be the G\"odel's number of \eqref{eq2.2-2}. Replacing the formula $y$ of the function $Sb(y, z, y)$ with the formula $n$, the formula \eqref{eq2.2-2} becomes

\begin{equation}\label{eq2.2-3}
\neg (\exists \ulcorner r \urcorner: \exists \ulcorner s \urcorner: (P(\ulcorner r \urcorner, \ulcorner s \urcorner) \land (\ulcorner s \urcorner = \ulcorner Sb(n,z, n) \urcorner ))
\end{equation}

In the presentation with the name of the formula and
the formula, we have

\begin{equation}\label{eq2.2-4}
\ulcorner G \urcorner \hspace{22pt} \neg (\exists \ulcorner r \urcorner: \exists \ulcorner s \urcorner: (P(\ulcorner r \urcorner, \ulcorner s \urcorner) \land (\ulcorner s \urcorner = \ulcorner Sb(n, z, n) \urcorner )))
\end{equation}
\begin{flushleft}
\hspace{77pt} {\scriptsize name of the formula} \hspace{60pt} {\scriptsize formula}
\end{flushleft}

%
%

\subsection{\textit{sub} and \textit{Sb}} \label{sec2.3}

\hspace{22pt}The steps of the procedure to obtain $sub(\ulcorner x \urcorner, \ulcorner z \urcorner, \ulcorner w \urcorner)$ are

\begin{itemize}
\item[(1)]  The entries are Gödel's numbers $\ulcorner x \urcorner$,  $\ulcorner z \urcorner$ e  $\ulcorner w \urcorner$.

\item[(2)] Locate the occurrences of $\ulcorner z \urcorner$ in the expanded Gödel's number of $\ulcorner x \urcorner$.

\item[(3)] Replace the occurrences of $\ulcorner z \urcorner$ with the expanded Gödel's number $\ulcorner w \urcorner$, generating the formula with the expanded Gödel's number $\ulcorner h \urcorner$.

\item[(4)] Calculate the Gödel's number of the expanded $\ulcorner h \urcorner$; this number is $sub(\ulcorner x \urcorner, \ulcorner z \urcorner, \ulcorner w \urcorner)$.
\end{itemize}

\textbf{An example:}

\hspace{22pt} In the following equality, the left-hand side is the Gödel number of $z = 0$ and the right-hand side is the expanded Gödel's number:
\begin{equation}\label{eq2.3-0}
497664000000 = 2^{17} \times 3^{5} \times 5^{6}
\end{equation}

\hspace{22pt} Let the set of consecutive primes be
\begin{equation}\label{eq2.3-1}
\{2, 3, 5, \ldots , 23, 29, \ldots, 107\}
\end{equation}
and the sets
\begin{equation}\label{eq2.3-2}
\{a, b, c, \ldots, p, q, \ldots, t\}
\end{equation}
\begin{equation}\label{eq2.3-3}
\{\alpha, \beta, \gamma\}
\end{equation}
where each element of \eqref{eq2.3-2} and \eqref{eq2.3-3} is one of the numbers assigned by Gödel to symbols ($\sim$, $\lor$, $\supset$, $\exists$, $=$, \ldots , $x$, $y$, $z$,  $\ldots$ ). Let the formula be $x$ whose Gödel's number is
\begin{equation}\label{eq2.3-4}
\ulcorner x \urcorner = 2^a \times 3^b \times 5^c \times \ldots \times 23^p \times 29^q \times \ldots \times 101^t
\end{equation}
Let the exponent $p$ in \eqref{eq2.3-4} be the symbol $z$ whose Gödel's number is 17. Assume that in \eqref{eq2.3-2} only $p$ is the free variable $z$. We have
\begin{equation}\label{eq2.3-5}
\ulcorner x \urcorner = 2^a \times 3^b \times 5^c \times \ldots \times 23^{17} \times 29^q \times \ldots \times 101^t
\end{equation}

\hspace{22pt} Let the formula $y$ have a Gödel's number that is
\begin{equation}\label{eq2.3-6}
\ulcorner w \urcorner = 2^{\alpha}\times 3^{\beta} \times 5^{\gamma}
\end{equation}
The result $sub(\ulcorner x \urcorner, \ulcorner z \urcorner, \ulcorner w \urcorner)$ is the number equal to
\begin{equation}\label{eq2.3-5}
2^a \times 3^b \times 5^c \times \ldots \times 23^{\alpha} \times 29^{\beta} \times 31^{\gamma} \times 37^{q} \times \ldots \times 107^t
\end{equation}

\hspace{22pt} The procedure for obtaining $Sb$ is different from the $sub$ procedure. Let the formulas be
\begin{equation}\label{eq2.3-7}
y = a b c \ldots p q \ldots t
\end{equation}
and
\begin{equation}\label{eq2.3-8}
w = \alpha \beta \gamma
\end{equation}
In \eqref{eq2.3-7} only $p$ is the free variable $z$ and there is only one type of free variable, that is,
\begin{equation}\label{eq2.3-9}
Sb(x,z,w) = a b c \ldots \alpha \beta \gamma q \ldots t
\end{equation}
The Gödel's number in the formula \eqref{eq2.3-9} is equal to \eqref{eq2.3-5}:
\begin{equation}\label{eq2.3-10}
\ulcorner Sb(x,z,w) \urcorner = 2^a \times 3^b \times 5^c \times \ldots \times 23^{\alpha} \times 29^{\beta} \times 31^{\gamma} \times 37^{q} \times \ldots \times 107^t
\end{equation}

%
%

\subsection{On the use of the \textit{sub} operator} \label{sec2.4}

\hspace{22pt}In section \ref{sec2.1}, the following hypothesis is made: “Assuming that the variable $\ulcorner y \urcorner$ of the operator $sub((\ulcorner y \urcorner, \ulcorner z \urcorner, \ulcorner y \urcorner)$ has the numerical value $\ulcorner n \urcorner$, the formula in \eqref{eq2.1-3} changes to $\ldots$”. This statement requires examination. Let's consider the formula $y$ given in \eqref{eq2.1-2}
\begin{equation}\label{eq2.4-1}
\ulcorner y\urcorner \hspace{44pt} \Delta z \Omega
\end{equation}
where $\Delta$ and $\Omega$ are sequences of symbols. Using the operator \textit{sub} we obtain

\begin{equation}\label{eq2.4-2}
sub(\ulcorner y \urcorner, \ulcorner z \urcorner, \ulcorner y \urcorner) \hspace{44pt} \Delta \Delta z \Omega \Omega
\end{equation}
\begin{flushleft}
\hspace{120pt} {\scriptsize name of the formula} \hspace{40pt}
{\scriptsize formula}
\end{flushleft}

In Nagel $\&$ Newman, in \eqref{eq2.1-3}, we have ``$\ldots$ [$sub(\ulcorner y \urcorner,17,\ulcorner y \urcorner)$] is still open-ended and indefinite, since it still contains the variable $y$. To make it definite, we need a numeral in place of a variable". Note that in \eqref{eq2.4-2} the formula also contains a free variable $z$.

%
%

\section{G\"odel's formula and the diagonal theorem} \label{sec3}

\begin{theorem} \label{teorema-1}
If the number of provable formulae is finite, then the formula $G$ doesn't have Gödel's number.
\end{theorem}
\textit{Proof:}

\hspace{22pt} Let be the set
\begin{equation}\label{eq2-9}
\{(\ulcorner r_1 \urcorner, \ulcorner s_1 \urcorner), (\ulcorner r_2 \urcorner, \ulcorner s_2 \urcorner), (\ulcorner r_3 \urcorner, \ulcorner s_3 \urcorner), \ldots, (\ulcorner r_k \urcorner, \ulcorner s_k \urcorner)\}
\end{equation}

where $\ulcorner s_i \urcorner$ is G\"odel's number of the proof of the formula of G\"odel's number $\ulcorner r_i \urcorner$. Let's consider that the number of provable
formulae $k$ is finite.

\hspace{22pt} Let's construct the formula $G$ as follows.

\begin{equation}\label{dois-10}
\neg (G = r_1) \land \neg (G = r_2) \land \ldots \land \neg (G = r_k) 
\end{equation}

\eqref{dois-10} in English means: \textit{The formula $G$ differs from all demonstrable formulas, therefore it is not demonstrable}.

\hspace{22pt} If $\ulcorner G \urcorner$ exists, then we can write
\begin{equation}\label{dois-11}
\ulcorner G \urcorner \hspace{22pt} \neg (G = r_1) \land \neg (G = r_2) \land \ldots \land \neg (G = r_k) 
\end{equation}
\begin{flushleft}
\hspace{77pt} {\scriptsize name of the formula} \hspace{60pt}
{\scriptsize formula}
\end{flushleft}
Examining \eqref{dois-11} it becomes evident that the Gödel's number of \eqref{dois-10} is not $\ulcorner G \urcorner$. The Gödel's number of the first sequence $\neg (G$ in the formula is already greater than $\ulcorner G \urcorner$.
\\
\hspace{10pt} $\Box$

\hspace{22pt} The diagonal theorem, which is the core of G\"odel's
proof, says ~\cite{Yanofsky}:

\begin{theorem} \label{teorema-2}
For any well-formed formula (wf) D(x) with x as its only free
variable, there exists a closed formula $\alpha$ such that
\begin{equation}\label{tres-1}
\vdash_{T} \ \alpha \longleftrightarrow D( \ulcorner \alpha
\urcorner )
\end{equation}
\end{theorem}

\begin{theorem} \label{teorema-3}
If there is an infinite set $\mathbb{S}$ of well formed formulae
with one free variable such that the elements of $\mathbb{S}$ have proofs in T, then the G\"odel's sentence does not have Gödel's number.
\end{theorem}

\hspace{22pt} The formula $D(x)$ used in G\"odel's proof is
\begin{equation}\label{tres-2}
D(s) = ( \forall r ) \neg P(r, s)
\end{equation}

Substituting \eqref{tres-2} in \eqref{tres-1} we obtain
\begin{equation}\label{tres-3}
\vdash _{T} \ G(x) \longleftrightarrow ( \forall r ) \neg P(r,
\ulcorner G(x) \urcorner)
\end{equation}

\hspace{22pt} Let
\begin{equation}\label{tres-4}
\mathbb{S} = \{ s_1, s_2, s_3, \ldots, s_k, \ldots \}
\end{equation}

be the set of all well formed formulas with one free variable that have proofs in $T$. The set $\mathbb{S}$ is not finite. We can construct the following sequence of symbols, denoted by
$\gamma$,
\begin{equation}\label{tres-5}
\neg ( Z =  s_1) \land \neg ( Z = s_2 ) \land \ldots \land \neg ( Z = s_k) \land \ldots
\end{equation}

The sequence $\gamma$ says with symbols of the list of symbols of $T$ that $Z$ is not provable. There are infinite
symbols in $\gamma$ and there is no way to express what $\gamma$ says with a number finite of symbols. Therefore the formula
\begin{equation}\label{tres-6}
( \forall r ) \neg P(r, \ulcorner G(x) \urcorner)
\end{equation}
has not a G\"odel's number (the number obtained is not finite). \hspace{10pt}
$\Box$

%
%

\section{Conclusion} \label{sec4}

\hspace{22pt} Based on the theorems in section \ref{sec3}, if there is a finite number of provable propositions (theorems) in the theory, then there is no Gödel's number for the formula $G$. If there is a non-finite number of provable propositions in the theory, then the Gödel's sentence does not have Gödel's number.


\begin{thebibliography}{99}
\addcontentsline{toc}{chapter}{Bibliography}
\bibitem{Tarski}  A. Tarsky, \textit{The Semantic Conception of Truth and the Foundations of Semantics}, Philosophy and Phenomenological Research \textbf{4}
(1944).
\bibitem{Hofstadter}  D. R. Hofstadter, \textit{G\"odel, Escher, Bach: an Eternal Golden Braid}, Vintage Books, (1980).
\bibitem{Nagel&Newman}  E. Nagel and J. R. Newman, \textit{Prova de G\"odel}, Editora Perspectiva, (1973).
\bibitem{Yanofsky}  Noson S. Yanofsky, \textit{A Universal Approach to Self-referential Paradoxes, Incompleteness and Fixed Points}, arXiv:math.LO/0305282 v1 19 May 2003.
\end{thebibliography}
\end{document}